\documentclass[12pt]{amsart}    
\usepackage{amscd}

\def\cal#1{\mathcal{#1}}

\def\NZQ{\Bbb}               
\def\NN{{\NZQ N}}

\def\ZZ{{\NZQ Z}}

\def\frk{\frak}               

\def\mm{{\frk m}}

\def\opn#1#2{\def#1{\operatorname{#2}}} 
\opn\chara{char}
\opn\length{\ell}
\opn\pd{pd}
\opn\rk{rk}
\opn\projdim{proj\,dim}
\opn\rank{rank}
\opn\depth{depth}
\opn\grade{grade}
\opn\height{height}
\opn\embdim{emb\,dim}
\opn\codim{codim}

\opn\Tr{Tr}
\opn\bigrank{big\,rank}
\opn\superheight{superheight}\opn\lcm{lcm}
\opn\trdeg{tr\,deg}%
\opn\reg{reg}
\opn\lreg{lreg}
\opn\div{div}
\opn\Div{Div}
\opn\cl{cl}
\opn\Cl{Cl}
\opn\Spec{Spec}
\opn\Supp{Supp}
\opn\supp{supp}
\opn\Sing{Sing}
\opn\Ass{Ass}
\opn\Ann{Ann}
\opn\Rad{Rad}
\opn\Soc{Soc}
\opn\Ker{Ker}
\opn\Coker{Coker}
\opn\Im{Im}
\opn\Hom{Hom}
\opn\Tor{Tor}
\opn\Ext{Ext}
\opn\End{End}
\opn\Aut{Aut}
\opn\id{id}

\opn\nat{nat}
\opn\pff{pf}
\opn\Pf{Pf}
\opn\GL{GL}
\opn\SL{SL}
\opn\mod{mod}
\opn\ord{ord}
\opn\aff{aff}
\opn\con{conv}
\opn\relint{relint}
\opn\st{st}
\opn\lk{lk}
\opn\cn{cn}
\opn\core{core}
\opn\vol{vol}
\opn\link{link}
\opn\star{star}
\opn\gr{gr}

\def\poly#1#2#3{#1[#2_1,\dots,#2_{#3}]}
\def\pot#1#2{#1[\kern-0.28ex[#2]\kern-0.28ex]}

\opn\dirlim{\underrightarrow{\lim}}
\opn\inivlim{\underleftarrow{\lim}}
\let\tensor=\otimes
\let\iso=\cong

\let\Dirsum=\bigoplus

\let\mcone= * 
\let\To=\longrightarrow
\def\Implies{\ifmmode\Longrightarrow \else
     \unskip${}\Longrightarrow{}$\ignorespaces\fi}
\def\implies{\ifmmode\Rightarrow \else
     \unskip${}\Rightarrow{}$\ignorespaces\fi}
\def\iff{\ifmmode\Longleftrightarrow \else
     \unskip${}\Longleftrightarrow{}$\ignorespaces\fi}

\let\:=\colon
\newtheorem{Theorem}{Theorem}[section]
\newtheorem{Lemma}{Lemma}[section]
\newtheorem{Corollary}{Corollary}[section]

\newtheorem{Remark}{Remark}

\newtheorem{Example}{Example}[section]

\let\epsilon\varepsilon
\let\kappa=\varkappa
\opn\inii{in}
\opn\inim{inm}
\opn\set{set}
\def\pnt{{\raise0.5mm\hbox{\large\bf.}}}

\begin{document}

\title{
A generalized Hochster's formula for
 local cohomologies of monomial ideals
}
\author{Yukihide Takayama}
\address{Yukihide Takayama, Department of Mathematical
Sciences, Ritsumeikan University, 
1-1-1 Nojihigashi, Kusatsu, Shiga 525-8577, Japan}
\email{takayama@se.ritsumei.ac.jp}
\date{September 12, 2004}

\def\Coh#1#2{H_{\mm}^{#1}(#2)}

\newcommand{\AppTh}{Theorem~\ref{approxtheorem} }
\def\da{\downarrow}
\newcommand{\ua}{\uparrow}
\newcommand{\namedto}[1]{\buildrel\mbox{$#1$}\over\rightarrow}
\newcommand{\bdel}{\bar\partial}
\newcommand{\proj}{{\rm proj.}}

\maketitle

\opn\Hilb{Hilb}

\begin{abstract}

The Hilbert series of local cohomologies for monomial ideals, which
are not necessarily square-free, is established.  As applications, we
give a sharp lower bound of the non-vanishing degree of local
cohomologies and also  a sharp lower bound of the positive integer 
$k$ of $k$-Buchsbaumness  for generalized Cohen-Macaulay monomial ideals.

\end{abstract}

\section*{Introduction}

Let $K$ be a field and 
let $S = \poly{K}{X}{n}$ be a polynomial ring with the standard
grading.
For a graded ideal $I\subset S$ we set $R = S/I$.  
We denote by $x_i$ the image
of $X_i$ in $R$ for $i=1,\ldots, n$ and set $\mm =(x_1,\ldots, x_n)$,
the unique graded maximal ideal.  Also $\Coh{i}{R}$ denotes the local
cohomology module of $R$ with regard $\mm$.

The aim of this paper is to show a generalization of Hochster's formula 
on local cohomologies for square-free monomial ideals (Stanley-Reisner
ideals) \cite{Hoc}
to monomial ideals that are not necessarily square-free.
The obtained formula, although its topological meaning is not clear
as compared to the original formula, tells much about the
non-vanishing degrees of the local cohomologies $\Coh{i}{R}$.  In
particular, we consider generalized Cohen-Macaulay monomial ideals.

A residue class ring $R$ 
is called {\em generalized Cohen-Macaulay} ring, 
or simply {\em generalized CM}, 
if $\Coh{i}{R}$ has finite length  for $i\ne \dim R$.
In this case, we will call the ideal $I\subset S$ a 
{\em generalized CM ideal}.
For a generalized CM ring $R$,
there exists an integer $k\in\ZZ$,
$k\geq 1$, such that $\mm^k\Coh{i}{R} = 0$ for $i\ne \dim R$. If this
condition holds, we will also call $R$, or $I\subset S$, {\em
$k$-Buchsbaum}. An ideal $I$ is generelized CM if and only if it is
$k$-Buchsbaum for some $k$. 
If $I$ is $k$-Buchsbaum but not
$(k-1)$-Buchsbaum, then we will call $I$ {\em strictly $k$-Buchsbaum}.
As a main application of our formula,
we will give a sharp bound of $k$ for the $k$-Buchsbaumness for a
generalized CM monomial ideals.

For a finite set $S$ we denote by $\mid S\mid$ the cardinarity of 
$S$, and, for sets $A$ and $B$, $A\subset B$ means that 
$A$ is a subset of $B$, which may be equal to $A$.

\thanks{The author thanks J\"urgen Herzog for 
valuable discussions and detailed comments on the early version of 
the paper.}

\section{Generalized Hochster's Formula}

We first consider an extension of Hochster's formula on local 
cohomologies of Stanley-Reisner ideals.

Let $I\subset S$  be a monomial ideal,
which is not necessarily square-free.
Then we have 
\begin{equation*}
	\Coh{i}{R}\iso H^i(C^\bullet)
\end{equation*}
where $C^\bullet$ is the $\check{C}$ech complex defined as follow:
\begin{equation*}
 C^\bullet : 
0\To C^0- \To C^1 \To \cdots \To C^n \To 0,
\qquad C^t = \Dirsum_{1\leq i_1<\cdots <i_t \leq n}R_{x_{i_1}\cdots x_{i_t}}.
\end{equation*}
and the differenctial $C^t\To C^{t+1}$ of this complex is induced by 
\begin{equation*}
(-1)^s nat : R_{x_{i_1}\cdots x_{i_t}} \To R_{x_{j_1}\cdots x_{j_{t+1}}}
\qquad \text{with } \{i_1,\ldots, i_t\} = \{j_1,\ldots, \hat{h}_s, \ldots, j_{t+1}\}
\end{equation*}
where $nat$ is the natural homomorphism to localized rings
and $R_{x_{i_1}\cdots x_{i_t}}$, for example, denotes localization by
$x_{i_1},\ldots, x_{i_t}$.

We can consider a $\ZZ^n$-grading to $\Coh{i}{R}$, $C^\bullet$
and $R_{x_{i_1}\cdots x_{i_t}}$ induced by the multi grading of 
$S$. See for example \cite{BH} for more detailed information about this complex.

Now we will consider the degree $a$ subcomplex $C^\bullet_a$ 
of $C^\bullet$ for any $a\in\ZZ^n$. Before that we will prepare the notation.
For a monomial ideal $I\subset S$, we denote by $G(I)$ the 
minimal set of monomial generators. Let $u = X_1^{a_1}\cdots X_n^{a_n}$ 
be a monimial with $a_i\geq 0$ for all $i$, 
then we define $\nu_j(u) = a_j$ for all $j=1,\ldots, n$
and $\supp(u)=\{i \mid a_i\ne 0\}$.
Now for $a\in\ZZ^n$, we set
$G_a = \{i \mid a_i < 0\}$ and $H_a = \{i \mid a_i > 0\}$.

\def\hochster#1{for all $u\in G(I)$ there exists $j\notin {#1}$
such that $\nu_j(u) > a_j \geq 0$}

\begin{Lemma}
\label{hochster:lemma1}
Let $x= x_{i_1}\cdots x_{i_r}$ with $i_1 < \cdots <i_r$ and 
set $F = \supp(x)$.
For all $a\in\ZZ^n$ we have $\dim_K (R_x)_a \leq 1$ and 
the following are equivalent
\begin{enumerate}
	\item [$(i)$] $(R_x)_a \iso K$ 
	\item [$(ii)$] 
	$F \supset G_a$ and \hochster{F}.
	\end{enumerate}
\end{Lemma}
Notice that the condition $a_i \geq 0$ 
in $(ii)$ is redundant because it follows from the condition $F\supset G_a$.
But it is written for the readers convenience.
\begin{proof} The proof of $\dim_K (R_x)_a \leq 1$ is  
verbatim the same as that of Lemma~5.3.6 (a) in \cite{BH}. 
Now we assume $(i)$, i.e.,
$(R_x)_a \ne 0$. This is equivalent to the condition that 
there exists a monomial $\sigma\in R$ and $\ell\in\NN$ such that 
\begin{enumerate}
\item [$(a)$] $x^m \sigma\ne 0$ for all $m\in\NN$, and 
\item [$(b)$] $\deg\displaystyle{\frac{\sigma}{x^\ell}} = a$,
\end{enumerate}
where $\deg$ denotes the multidegree.
We know from $(b)$ that we have $F \supset G_a$ 
because a negative degree $a_i (<0)$ in $a$ must 
come from the denominator of the fraction $\sigma/x^\ell$ and 
$F = \supp(x^\ell)$.
Now we know that $(a)$ is equivalent to the following 
condition:
for all $u\in G(I)$ and for all $m\in\NN$ we have
$u\not| (X_{i_1}^m\cdots X_{i_r}^m)(X_1^{b_1}\cdots X_n^{b_n})$
where we set $\sigma = x_1^{b_1}\cdots x_n^{b_n}$ with some integers $b_j\geq 0$,
$j=1,\ldots, n$. This is 
equivalent to the following:
for all $u\in G(I)$ there exists $i\notin F$
such that $\nu_i(u) > b_i$. Furthermore, 
we know from the condition $F \supset G_a$ 
that we have $a_i = b_i$ for $i\notin F$ since
by $(b)$ non-negative degrees in $a$ must come from $\sigma$.
Consequently we obtain $(ii)$.

Now we show the converse. Assume that we have $(ii)$. 
Set $\tau = \prod_{i\in H_a}x_i^{a_i}$   
and $\rho = \prod_{i\in G_a}x_i^{-a_i}$. 
Then since $F \supset G_a$ there exists $\ell \in \NN$ and 
a monomial $\sigma$ (in $R$) 
such that 
\begin{equation}
\label{hochster:lemma1-eq1}
	x^\ell = \rho \sigma   
\end{equation}
We show that 
$\displaystyle{\frac{\sigma\tau}{x^\ell}}\ne 0$ 
in $R_x$.
$\displaystyle{\frac{\sigma\tau}{x^\ell}}\ne 0$ 
is equivalent to the condition
that 
$x^m(\sigma\tau)\ne 0$ 
for all $m\in\NN$, as in the above discussion. 
This is equivalent to the condition
\begin{equation}
\label{hochster:lemma1-eq2}
	\text{for all } u \in G(I)\; \text{ there exists } i\notin F
\mbox{ such that } \nu_i(u) > b_i
\end{equation}
where we set 
$\sigma\tau = x_1^{b_1}\cdots x_n^{b_n}$ 
for some integers $b_j\geq 0$, $j=1,\ldots, n$.
But by (\ref{hochster:lemma1-eq1}) we have 
$i\notin\supp(\sigma)$ 
for
$i\notin F$, so that 
$b_i = \nu_i(\tau) = a_i(>0)$ 
(i.e., 
$i\in H_a$)
or $a_i = b_i = 0$ (i.e., $i\notin H_a\cup G_a$). Hence
(\ref{hochster:lemma1-eq2}) is exactly the condition
that \hochster{F},  which is assured by the assumption.
Thus we have 
$\displaystyle{\frac{\sigma\tau}{x^\ell}}\ne 0$ in $R_x$.
Therefore
\begin{equation*}
   \deg\displaystyle{\frac{\sigma\tau}{x^\ell}}
=  \deg\displaystyle{\frac{\sigma\tau}{\rho\sigma}}
=  \deg\prod_{i\in H_a\cup G_a}x_i^{a_i}
= \deg x^a = a
\end{equation*}
as required.
\end{proof}

Let $a\in\ZZ^n$. By Lemma~\ref{hochster:lemma1} we see that 
$(C^i)_a$ has a $K$-linear basis
\begin{equation*}
   \{
	b_F 
	: F\supset G_a, \text{ and } \mbox{\hochster{F}}
   \}.
\end{equation*}
Restricting the differentation of $C^\bullet$ to the
$a$th graded piece we obtain a complex $(C^\bullet)_a$
of finite dimensional $K$-vector spaces with differentation
$\partial : (C^i)_a\To (C^{i+1})_a$ given
by $\partial(b_F) = \sum(-1)^{\sigma(F, F')}b_{F'}$
where the sum is taken over all $F'$ such that 
$F' \supset F$ with $\vert F'\vert = i+1$
and  \hochster{F'}. Also we define $\sigma(F, F') =s$
if  $F'= \{j_0,\ldots, j_i\}$ and 
$F = \{j_0,\ldots, \hat{j}_s,\ldots, j_i\}$.
Then we describe the $a$th component of the local cohomology
in terms of this subcomplex:
$\Coh{i}{R}_a \iso H^i(C^\bullet)_a = H^i(C^\bullet_a)$.

Now we fix our notation on simplicial complex. A simplicial 
complex $\Delta$ on a finite set $[n] = \{1,\ldots, n\}$
is a collection of subsets of $[n]$ such that $F\in\Delta$
whenever $F\subset G$ for some $G\in \Delta$.
Notice that we do not assume the condition
that $\{i\}\in \Delta$ for $i=1,\ldots, n$. 
We define $\dim F = i$ if $\mid F\mid = i+1$
and $\dim \Delta = \max\{\dim F \mid F\in\Delta\}$,
which will be called the dimension of $\Delta$.
If we assume a 
linear order on $[n]$, say $1<2<\cdots<n$, then we will call
$\Delta$ {\em oriented}, and in this case we always denote
an element $F=\{i_1,\ldots, i_k\}\in\Delta$ 
with the orderd sequence $i_1 <\ldots < i_k$.
For a given oriented simplicial complex
of dimension $d-1$, 
we denote by ${\cal C}(\Delta)$ the augumented oriented chain complex 
of $\Delta$:
\begin{equation*}
{\cal C}(\Delta) :
0\To {\cal C}_{d-1}\overset{\partial}{\To}
     {\cal C}_{d-2}\To\cdots\To
     {\cal C}_0\overset{\partial}{\To}
     {\cal C}_{-1}\To 0
\end{equation*}
where 
\begin{equation*}
   {\cal C}_i = \Dirsum_{F\in\Delta, \dim F=i}\ZZ F
\qquad\text{and}\qquad
   \partial F = \sum_{j=0}^i(-1)^jF_j
\end{equation*}
for all $F\in \Delta$. Here we define
$F_j = \{i_0,\ldots, \hat{i}_j,\ldots, i_k\}$
for $F = \{i_0,\ldots, i_k\}$.
Now for an abelian group $G$, we define
the $i$th reduced simplicial homology 
$\tilde{H}_i(\Delta; G)$ 
of $\Delta$
to be the $i$th homology 
of the complex ${\cal C}(\Delta)\tensor G$ for 
all $i$. Also we define the 
$i$th reduced simplicial cohomology of $\Delta$
$\tilde{H}^i(\Delta; G)$ to be the $i$th cohomology
of the dual chain complex $\Hom_\ZZ({\cal C}(\Delta), G)$
for all $i$.
Notice that we have 
\begin{equation*}
\tilde{H}_{-1}(\Delta; G) 
 = \left\{
	\begin{array}{ll}
	   G  &  \text{if $\Delta = \{\emptyset\}$} \\
	   0  &  \text{otherwise}
	\end{array}
   \right.,
\end{equation*}
and if $\Delta = \emptyset$ then $\dim \Delta = -1$ and 
$\tilde{H}_{i}(\Delta; G)=0$ for all $i$.

Now we will establish an isomorphism between the complex
$(C^\bullet)_a$, $a\in\ZZ^n$, and a dual chain complex.
For any $a\in\ZZ^n$, we define a simplicial complex
\begin{equation*}
\Delta_a = 
\left\{ F- G_a
	\;\vert\;
	\begin{array}{l}
	 F \supset G_a \mbox{ and }\\
        \mbox{\hochster{F}}
	\end{array}
\right\}.
\end{equation*}
Notice that we may have $\Delta_a=\emptyset$ for some $a\in\ZZ^n$.

\begin{Lemma}
\label{hochster:lemma2}
For all $a\in\ZZ^n$ there exists an isomorphism
of complexes
\begin{equation*}
  \alpha^\bullet 
:  (C^\bullet)_a \To 
   \Hom_\ZZ({\cal C}(\Delta_a)[-j-1], K)
\qquad j = \vert G_a\vert
\end{equation*}
where ${\cal C}(\Delta_a)[-j-1]$ means shifting the 
homological degree 
of ${\cal C}(\Delta_a)$ by $-j-1$.
\end{Lemma}
\begin{proof}
The assignment $F\mapsto F - G_a$ induces an isomorphism
$\alpha^\bullet : (C^\bullet)_a \To \Hom_\ZZ({\cal C}(\Delta_a)[-j-1], K)$ 
of $K$-vector spaces such that $b_F \mapsto \varphi_{F-G_a}$,
where
\begin{equation*}
   \varphi_{F'}(F^")
= \left\{
	\begin{array}{ll}
	    1 & \mbox{if $F' = F^"$} \\
            0 & \mbox{otherwise.}
	\end{array}
  \right. 
\end{equation*}
That this is a homomorphism of complexes can be checked 
in a straightforward way.
\end{proof}

Now we come to our main theorem.

\begin{Theorem}
\label{hochster:theoremTheFormula}
Let $I\subset S = \poly{K}{X}{n}$ be a monomial ideal. Then
the multigraded Hilbert series of the local cohomology modules of $R = S/I$ 
with respect to the $\ZZ^n$-grading is given by
\begin{equation*}
\Hilb(\Coh{i}{R}, {\bf t})
=\sum_{F\in\Delta}
  \sum
   \dim_K\tilde{H}_{i-\vert F\vert -1}(\Delta_a; K) {\bf t}^a
\end{equation*}
where 
${\bf t} = t_1\cdots t_n$,
the second sum runs over 
$a\in\ZZ^n$ such that 
$G_a = F$ and $a_j\leq \rho_j-1$, $j=1,\ldots,n$,
with 
$\rho_j = \max\{\nu_j(u)\;\vert\; u\in G(I)\}$
for $j=1,\ldots, n$, and $\Delta$ is the 
simplicial complex corresponding to the Stanley-Reisner
ideal $\sqrt{I}$.
\end{Theorem}

\begin{proof}
By Lemma~\ref{hochster:lemma2} and universal coefficient theorem
for simplicial (co)homology, we have 
\begin{eqnarray*}
\Hilb(\Coh{i}{R},{\bf t})
& = & \sum_{a\in\ZZ^n}\dim_K\Coh{i}{R}_a {\bf t}^a \\
& = & \sum_{a\in\ZZ^n}\dim_K\tilde{H}_{i-\vert G_a\vert -1}(\Delta_a;K) 
                                         {\bf t}^a. \\
\end{eqnarray*}
It is clear from the definition that $\Delta_a= \emptyset$
if for all $j\notin G_a$ we have $a_j \geq \rho_j$.
In this case, we have 
$\dim_K\tilde{H}_{i-\vert G_a\vert -1}(\Delta_a;K) =0$.
Thus we obtain
\begin{eqnarray*}
\Hilb(\Coh{i}{R},{\bf t})
& = & \sum_{{\tiny
		\begin{array}{c}
	              a\in\ZZ^n\\
		   a_j\leq \rho_j-1\\
	           j=1,\ldots, n
		\end{array}
             }}
        \dim_K\tilde{H}_{i-\vert G_a\vert -1}(\Delta_a;K) 
                                         {\bf t}^a. \\
\end{eqnarray*}
Now if $\Delta_a\ne \emptyset$, we must have $(G_a - G_a =)\;
\emptyset\in \Delta_a$, i.e., $G_a$ must be a subset of 
$\{1,\ldots, n\}$ such that \hochster{G_a}, and this condition is
equivalent to "$G_a \not\supset \supp(u)$ for all $u\in G(I)$",
which can further be refined as "$G_a$ is not 
a non-face of $\Delta$, i.e., $G_a\in \Delta$".
Thus we finally obtain the required formula.
\end{proof}

The original Hochster's formula is a special case of 
Theorem~\ref{hochster:theoremTheFormula}. 
For a simplicial complex $\Gamma$ and $F\in \Gamma$,
we define $\lk_\Gamma F = \{G \vert F\cup G\in\Gamma, F\cap G=\emptyset\}$
and $\st_\Gamma F = \{G \vert F\cup G\in\Gamma\}$.

\begin{Corollary}[Hochster]
Let $\Delta$ be a simplicial complex and 
let $K[\Delta]$ be the Stanley-Reisner ring 
corresponding to $\Delta$. Then we have
\begin{equation*}
\Hilb(\Coh{i}{K[\Delta]}, {\bf t})
=\sum_{F\in\Delta}
   \dim_K\tilde{H}_{i-\vert F\vert -1}(\lk_\Delta F; K) 
	\prod_{j\in F}
	   \displaystyle{\frac{t_j^{-1}}{1- t_j^{-1}}}.
\end{equation*}
\end{Corollary}
\begin{proof}
By Theorem~\ref{hochster:theoremTheFormula} we have 
\begin{equation*}
\Hilb(\Coh{i}{R}, {\bf t})
=\sum_{F\in\Delta}
 \sum_{{\tiny \begin{array}{c}
	a\in\ZZ^n_{-} \\
        G_a = F
	\end{array}}}
   \dim_K\tilde{H}_{i-\vert F\vert -1}(\Delta_a; K) {\bf t}^a
\end{equation*}
where $\ZZ^n_{-} = \{a\in\ZZ^n \vert a_j\leq 0\text{ for }j=1,\ldots,n\}$
and 
\begin{eqnarray*}
\Delta_a &=& 
\left\{ F- G_a
	\;\vert\;
	\begin{array}{l}
	 F \supset G_a, \mbox{ and for all $u\in G(I)$ there exists $j\notin F$}\\
        \mbox{
           such that $j\in \supp(u)$ and $j\notin H_a\cup G_a$
             }
	\end{array}
\right\}.\\
&=& 
\left\{ F- G_a
	\;\vert\;
	 F \supset G_a, \text{ and for all $u\in G(I)$ we have
          $H_a\cup F\not\subset \supp(u)$}
\right\}.\\
&=& 
\left\{ L
	\;\vert\;
         L\cap G_a=\emptyset, 
	 L\cup G_a\cup H_a\in\Delta
\right\} = \lk_{\st_\Delta H_a}G_a.
\end{eqnarray*}
Then the rest of the proof is exactly as in Theorem~5.3.8 \cite{BH}.
\end{proof}

\section{Applications}

In this section, we give some application of 
Theorem~\ref{hochster:theoremTheFormula}. 
We define $a_i(R) = \max\{j \vert \Coh{i}{R}_j\ne 0\}$
if $\Coh{i}{R}\ne 0$ and $a_i(R) = -\infty$ if $\Coh{i}{R}=0$.
Similarly, we define 
and $b_i(R) = \inf\{j \vert \Coh{i}{R}_j\ne 0\}$
if $\Coh{i}{R}\ne 0$ and  $b_i(R) = +\infty$ if 
$\Coh{i}{R} =0$.

Recall that 
$\rho_j = \max\{\nu_j(u)\;\vert\; u\in G(I)\}$
for $j=1,\ldots, n$.

\begin{Corollary}
\label{hochster:corollary1}
Let $I\subset S = \poly{K}{X}{n}$ be a monomial ideal.
Then $a_i(R) \leq \sum_{j=1}^n \rho_j - n$ for all $i$.
\end{Corollary}
\begin{proof}
By Theorem~\ref{hochster:theoremTheFormula},
the terms in 
$\Hilb(\Coh{i}{R}, {\bf t})$ 
with the highest total degree are 
at most $\dim_K\tilde{H}_{i-\vert F\vert -1}(\Delta_a;K){\bf t}^a$ 
with $a_j = \rho_j-1$ for all $j$. Thus the total degree is 
at most $\sum_{j}\rho_j - n$.
\end{proof}

From Corollary~\ref{hochster:corollary1}, we can recover the 
following well known result.

\begin{Corollary}
Let $I\subset S$ be a Stanley-Reisner ideal. Then 
$a_i(R) \leq 0$ for all $i$.
\end{Corollary}
\begin{proof}
If $I$ is square-free, then $\rho_j\leq 1$ for $j=1,\ldots, n$.
\end{proof}

For a Stanley-Reisner generalized Cohen-Macaulay ideal $I\subset S$
with $\dim R =d$,
it is well known that it is Buchsbaum and $b_i(R)\geq 0$ for all $i(\ne d)$. 
The following theorem extends this result to monomial ideals in general.

\begin{Theorem}
\label{hochster:theoremGenCMCase}
Let $I\subset S=\poly{K}{X}{n}$ be a generalized CM
monomial ideal.
Then $b_i(R) \geq 0$ for all $i < \dim R$.
\end{Theorem}
\begin{proof}
Let $d=\dim R$.
Assume that there exists $i$ and $j$ with $0\leq i<d$ and $j<0$
such that $\Coh{i}{R}_j\ne 0$. Then
by Theorem~\ref{hochster:theoremTheFormula} there exists
$a\in \ZZ^n$ such that 
\begin{enumerate}
\item [$(i)$] $\sum_{k=1}^n a_k = j < 0$, in particular $G_a\ne \emptyset$,
and
\item [$(ii)$] $\dim_K\tilde{H}_{i-\vert G_a\vert -1}(\Delta_a;K)\ne 0$,
in particular $\Delta_a \ne \emptyset$.
\end{enumerate}
Now observe that by the definition of $\Delta_a$, the coditions 
$(ii)$ is  independent of the values of $a_j$ for $j\in G_a$.
This means that the total degree $j = \sum_{k=1}^n a_k$ can be 
any negative integer so that $\Coh{i}{R}$ is not of finite length,
which contradicts the assumption.
\end{proof}

\begin{Remark}
{\em 
By Kodaira Vanishing Theorem $($Corollary~2.4 \cite{HK}$)$,
we have $b_i(R)\geq 0$ for $i\ne \dim R$ if $\chara(K)=0$ and 
$R$ is a normal domain and has an
isolated singularity at $\mm$.
Theorem~\ref{hochster:theoremGenCMCase} is a case that is 
not covered by Kodaira Vanishing Theorem.}
\end{Remark}

For a generalized Cohen-Macaulay ideal $I\subset S$ with $d = \dim R$,
there exists some positive integer $k$ such that $\mm^k \Coh{i}{R}=0$
for all $i\ne d$. Then we refer  $R$ as a $k$-Buchsbaum ring. 
Now we consider the question: what is the lower bound of $k$?

\begin{Theorem}
\label{hochster:theoremTheBound}
Let $I\subset S=\poly{K}{X}{n}$ be a generalized CM
monomial ideal. Then $R=S/I$ is 
$\left(\sum_{j=1}^n \rho_j -n+1\right)$-Buchsbaum.
\end{Theorem}
\begin{proof}
$R$ is $\max_{i\ne d}(a_i(R) - b_i(R) + 1)$-Buchsbaum.
The required result follows
immediately from 
Corollary~\ref{hochster:corollary1} and 
Theorem~\ref{hochster:theoremGenCMCase}.
\end{proof}

We can immediately recover the following well known result.

\begin{Corollary}
\label{hochster:SRis1Bbm}
Let $I\subset S$ be a Stanley-Reisner ideal. 
If $R$ is generalized Cohen-Macaulay, then $R$ is $1$-Buchsbaum.
\end{Corollary}

In fact, it is well-known that 
a generalized Cohen-Macaulay Stanley-Reisner ideal
is Buchsbaum, which is stronger than $1$-Buchsbaumness.

The bound of $k$-Buchsbaumness given in
Theorem~\ref{hochster:theoremTheBound} is best possible.
In fact, we can construct strictly $(\sum_{j=1}^n \rho_j -n+1)$-Buchsbaum
ideals  as in the following example.

\begin{Example}
{\em 
Let $I\subset S$ be a Stanley-Reisner Buchsbaum ideal.
Notice that such ideals can be constructed with the method
presented in \cite{BjoHib} and $\Coh{i}{S/I}$ ($i\ne \dim R$)
is a $K$-vector space for $i\ne \dim R$.

Now consider a $K$-homomorphism
\begin{equation*}
  \varphi : S \To S, \qquad X_i \longmapsto X_i^{a_i}\; (i=1,\ldots, n)
\end{equation*}
where $a = (a_1,\ldots, a_n) \in\ZZ^n$ with $a_i\geq 1$ for $i=1,\ldots, n$.
We define $\varphi(M) = M\tensor_S \hbox{}^\varphi S$ for a $S$-module $M$,
where a left-right $S$-module $\hbox{}^\varphi S$ is equal to $S$ as a set,
it is a right $S$-module in the ordinary sense and 
its left $S$-module structure is determined by $\varphi$.
Then we have 
\begin{enumerate}
\item $\varphi(S/I) = S/\varphi(I)S$,
\item $\varphi$ is an exact functor.
\end{enumerate}
Thus, for $i\ne \dim R$, 
we have $\Coh{i}{S/\varphi(I)S}\iso \varphi(\Coh{i}{S/I})$ and 
since $\Coh{i}{S/I}$ is a direct sum of $S/\mm$, 
$\Coh{i}{S/\varphi(I)S}$ is a direct sum of 
$S/(X_1^{a_1},\ldots, X_n^{a_n})$. Then we know that 
$\mm^k\Coh{i}{S/I} = 0$ but  $\mm^{k-1}\Coh{i}{S/I}\ne 0$ 
with $k = \sum_{j=1}^n \rho_j - n +1 = \sum_{j=1}^n a_j - n +1$.
}
\end{Example}

\begin{Remark}
{\em 
Bresinsky and Hoa gave a bound for $k$-Buchsbaumness
for ideals generated by monomials and binomials 
$($Theorem~4.5 \cite{BreHoa}$)$.
For monomial ideals our bound is stronger than
that of Bresinsky and Hoa.}
\end{Remark}

Finally, we consider vanishing cohomological dimensions of generalized
CM monomial ideals.
Recall that Castelnuovo-Mumford regularity of the ring $R$ 
is defined by 
\begin{equation*}
   \reg(R) = \max\{i+j \vert \Coh{i}{R}_j\ne 0\}.
\end{equation*}
Let $r = \reg(R)$. Then we have $\Coh{i}{R}_j =0$ 
for $j > r -i$. 
Then we have

\begin{Corollary}
Let $I\subset S$ be a generalized CM monomial ideal
with $d =\dim R$ and  $r = \reg(R)$. 
Then $\Coh{i}{R} =0$ for $r+1 \leq i<d$.
In particular, if $I$ has a $q$-linear resolution, 
we have  $\Coh{i}{R} =0$ for $q \leq i<d$.
\end{Corollary}
\begin{proof}
First part is clear from Theorem~\ref{hochster:theoremGenCMCase}.
If $R$ has a $q$-linear resolution, we have $\reg(R) = q-1$.
Thus the second statement also follows immediately.
\end{proof}


\end{document}